# НЕРАВЕНСТВА ТИПА БЕРНШТЕЙНА В ПОДПРОСТРАНСТВАХ ИНВАРИАНТНЫХ ОТНОСИТЕЛЬНА СДВИГА

В.Ф. Бабенко*'**, А.А. Лигун, С.А. Спектор*.

Украина. *Днепропетровский национальный университет

** Институт прикладной математики и механики НАН Украины

Решение многих задач математического анализа на классах периодических функций базируется на неравенствах для тригонометрических многочленов. Причем особое значение имеют неравенства, дающие оценку нормы производной тригонометрического многочлена через норму самого многочлена. Первое неравенство такого рода было получено в 1912 году С.Н.Бернштейном.

Классическое неравенство Бернштейна [1,122]

$$\| \tau_n' \|_C \le n \| \tau_n \|_C, \quad k = 1, 2, \ldots$$

для тригонометрического полинома $\tau(t)$ порядка n, имеет важные приложения в теории аппроксимации и других областях математики.

Точные неравенства типа Бернштейна для любого полинома $\tau(t) \in \mathrm{T}_{2n+1}$ известны в следующих случаях.

Неравенство

$$\left\| \tau_n^{(k)} \right\|_1 \le n^k \frac{\left\| \cos(\square) \right\|_1}{\left\| \cos(\square) \right\|_p} \left\| \tau_n \right\|_p \quad , \ \forall n, k = 1, 2 \ldots \ и \ 1 \le p \le \infty ,$$

установлено А.А.Лигуном [2].

Неравенства

$$\sqrt{\frac{2\pi}{2n+1}} \, \| \tau \|_C \le \| \tau \|_2 , \quad \forall \, n \in \mathrm{N}$$

$$\| \tau \|_p \le \left( \frac{n q_0 + 1}{2\pi} \right)^{1/q - 1/p} \| \tau \|_q , \ (где \ 1 \le q \le p \le \infty , \ а \ q_0 \ наименьшее \ четное \ число \ большее \ или \ равное \ q) ,$$

доказаны С.М.Никольским [3]. Также, более детальный обзор точных неравенств типа Бернштейна можно найти, например, в 3 главе работы [1].

Мы будем рассматривать вопрос о неравенствах типа Бернштейна для подпространств пространства $L_2(R)$ и $L_2(R^n)$, инвариантных относительно сдвига.

Пусть задана функция $\varphi(x) \in L_2(R)$. Рассмотрим ортонормированную систему функций

$$\{\varphi_h(x)\} = \{\varphi(\bullet - h) : h \in Z\}.$$

Очевидно, $\sum_{\ell \in R} \left|\hat{\varphi}(\omega + 2\pi\ell)\right|^2 = 1$ [4, 266-271].

Пусть $W$ - подпространство пространства $L_2(R)$ порожденное сдвигами функции $\varphi(x)$. Любую функцию $f \in W$ можно представить в виде суммы сходящегося в $L_2(R)$ ряда

$$f = \sum_{\gamma \in Z} c_\gamma \varphi(\bullet - h).$$

Преобразование Фурье функции $f$ зададим в виде

$$f(\omega) = \frac{1}{\sqrt{2\pi}} \int_R f(x) e^{-i\omega x} dx \quad .$$

Тогда,

$$\hat{f} = \sum_{\gamma \in Z} c_\gamma \hat{\varphi}(\omega) e^{-i\gamma\omega} = m_f(\omega) \hat{\varphi}(\omega), \text{ где } m_f(\omega) = \sum_{\gamma \in Z} c_\gamma e^{-i\gamma\omega} .$$

Нами доказана следующая

**Теорема 1.** Зададим функцию $\varphi(x)$ такую, что

$$\sup_{\omega} \sum_{\ell \in R} |\omega + 2\pi\ell|^{2k} \left|\hat{\varphi}(\omega + 2\pi\ell)\right|^2 < \infty .$$

Тогда, для любой функции $f \in W$, $\forall k \in Z$ выполняется точное неравенство

$$\| f^{(k)} \|_2^2 \le \sup_{\omega} \sum_{\ell \in R} |\omega + 2\pi\ell|^{2k} \left|\hat{\varphi}(\omega + 2\pi\ell)\right|^2 \| f \|_2^2 \qquad (1)$$

**Доказательство.** Воспользовавшись а равенство Парсеваля:

$$\int_R | f(x)|^2 dx = \frac{1}{2\pi} \int_R | \hat{f}(\omega)|^2 d\omega,$$

получим

$$\| f^{(k)}\|_2^2=\frac{1}{2\pi}\left\|(\mathrm{i}\omega)^k\hat f(\omega)\right\|_2^2=\frac{1}{2\pi}\int_R|i\omega|^{2k}\left|\hat f(\omega)\right|^2 d\omega=$$

$$=\frac{1}{2\pi}\int_R\omega^{2k}\left|\sum_{\gamma\in Z}c_\gamma e^{-\mathrm{i}\gamma\omega}\hat\varphi(\omega)\right|^2 d\omega=$$

$$=\frac{1}{2\pi}\sum_{\ell\in R}\int_0^{2\pi}|\omega+2\pi\ell|^{2k}\left|\sum_{\gamma\in Z}c_\gamma e^{-\mathrm{i}\gamma\omega}\hat\varphi(\omega+2\pi\ell)\right|d\omega=$$

$$=\frac{1}{2\pi}\int_0^{2\pi}\sum_{\ell\in Z}|\omega+2\pi\ell|^{2k}\left|\sum_{\gamma\in Z}c_\gamma e^{-\mathrm{i}\gamma\omega}\hat\varphi(\omega+2\pi\ell)\right|^2 d\omega=$$

$$=\frac{1}{2\pi}\int_0^{2\pi}|m_f(\omega)|^2\sum_{\ell\in Z}|\omega+2\pi\ell|^{2k}\left|\hat\varphi(\omega+2\pi\ell)\right|^2 d\omega\,.$$

Заметим, что

$$\| f\|_2^2=\frac{1}{2\pi}\int_R\left|\hat f(\omega)\right|^2 d\omega=\frac{1}{2\pi}\int_0^{2\pi}\sum_{\ell\in R}\left|\hat\varphi(\omega+2\pi\ell)\right|^2\left|\sum_{\gamma\in Z}c_\gamma e^{-\mathrm{i}\gamma\omega}\right|^2 d\omega=$$

$$=\frac{1}{2\pi}\int_0^{2\pi}|m_f(\omega)|^2\sum_{\ell\in R}\left|\hat\varphi(\omega+2\pi\ell)\right|^2=\frac{1}{2\pi}\int_0^{2\pi}|m_f(\omega)|^2 d\omega.$$

Таким образом, получим

$$\| f^{(k)}\|_2^2=\frac{1}{2\pi}\int_0^{2\pi}|m_f(\omega)|^2\sum_{\ell\in Z}|\omega+2\pi\ell|^{2k}\left|\hat\varphi(\omega+2\pi\ell)\right|^2 d\omega\le$$

$$\le\sup_\omega\sum_{\ell\in\mathrm{R}}|\omega+2\pi\ell|^{2k}\left|\hat\varphi(\omega+2\pi\ell)\right|^2\| f\|_2^2\,.$$

Докажем теперь, что константа в правой части неравенства (1) неулучшаема. Для этого нужно выбрать такое семейство функций $\mathrm{f_n}, n\in\mathrm{Z}$, чтобы

$$\frac{\| f_n^{(k)}\|_2^2}{\| f_n\|_2^2}\xrightarrow[n\to\infty]{}\sup_\omega\frac{\sum_{\ell\in R}|\omega+2\pi\ell|^{2k}\left|\hat\varphi(\omega+2\pi\ell)\right|^2}{\sum_{\ell\in R}\left|\hat\varphi(\omega+2\pi\ell)\right|^2}\,.$$

Пусть

$$\Phi_\mathrm{n}(\omega)=\frac{\sin^2\left(\mathrm{n}+\frac{1}{2}\right)\omega}{(\mathrm{n}+1)\sin^2\frac{1}{2}\omega}\quad -$$

ядро Фейера.

Положим

$$m_n(\omega)=\sqrt{\Phi_n(\omega)} \ .$$

Определим последовательность $f_{\text{n}}$ следующим образом:

$$\hat{f}_n(\omega)=m_n(\omega)\hat{\varphi}(\omega)$$

Заметим, что

$$\| f \|_2^2=\frac{1}{2\pi}\int_0^{2\pi}\left|m_f(\omega)\right|^2 d\omega=\frac{1}{2\pi}\int_0^{2\pi}\left|\Phi_{\text{n}}(\omega)\right|d\omega=1 \qquad (2)$$

Рассмотрим

$$\frac{\| f_n^{(k)}\|_2^2}{\| f_n\|_2^2}=\frac{\left\|(\mathrm{i}\omega)^k \hat{f}_n(\omega)\right\|_2^2}{\left\|\hat{f}_n\right\|_2^2}=\frac{\left\|(i\omega)^k\sqrt{\Phi_n(\omega)}\,\hat{\varphi}(\omega)\right\|_2^2}{\left\|\sqrt{\Phi_n(\omega)}\,\hat{\varphi}(\omega)\right\|_2^2}=$$

$$=\frac{\int_{\mathrm{R}}|\omega|^{2k}\Phi_{\mathrm{n}}(\omega)|\hat{\varphi}(\omega)|^2 d\omega}{\int_{\mathrm{R}}\Phi_{\mathrm{n}}(\omega)|\hat{\varphi}(\omega)|^2 d\omega}=\frac{\sum_{\ell\in R}\int_0^{2\pi}|\omega+2\pi\ell|^{2k}\Phi_n(\omega)\left|\hat{\varphi}(\omega+2\pi\ell)\right|^2 d\omega}{\sum_{\ell\in R}\int_0^{2\pi}\Phi_n(\omega)\left|\hat{\varphi}(\omega+2\pi\ell)\right|^2 d\omega}=$$

$$=\frac{\int_0^{2\pi}\sum_{\ell\in R}|\omega+2\pi\ell|^{2k}\Phi_n(\omega)\left|\hat{\varphi}(\omega+2\pi\ell)\right|^2 d\omega}{\int_0^{2\pi}\sum_{\ell\in R}\Phi_n(\omega)\left|\hat{\varphi}(\omega+2\pi\ell)\right|^2 d\omega}\xrightarrow[n\to\infty]{}$$

$$\xrightarrow[\mathrm{n}\to\infty]{}\sup_{\omega}\frac{\sum_{\ell\in\mathrm{R}}|\omega+2\pi\ell|^{2\mathrm{k}}|\hat{\varphi}(\omega+2\pi\ell)|^2}{\sum_{\ell\in\mathrm{R}}|\hat{\varphi}(\omega+2\pi\ell)|^2}=$$

$$=\sup_{\omega}\sum_{\ell\in\mathrm{R}}|\omega+2\pi\ell|^{2k}\left|\hat{\varphi}(\omega+2\pi\ell)\right|^2 \ .$$

Теорема доказана.

Рассмотрим систему функций нецелочисленных сдвигов

$$\{\varphi_h(x)\}=\left\{\varphi\left(\frac{\bullet-h}{\mathrm{a}}\right): h,a\in \mathrm{Z}\right\}.$$

Соответственно теперь

$$W=\left\{ f\in L_2(R): f=\sum_{\gamma\in Z} c_\gamma \varphi\left(\frac{x-\gamma}{a}\right)\right\} \text{ -}$$

подпространство пространства $L_2(R)$ порожденное нецелочисленными сдвигами функции $\varphi(x)$.

Для любой функции $f\in W$ имеет место следующая

**Теорема 2.** Зададим функцию $\varphi(x)$ такую, что $\forall\omega\in R$,

$$\sup_{\omega}\sum_{\ell\in R}\left|\frac{\omega+2\pi\ell}{a}\right|^{2k}\left|\hat{\varphi}\left(\frac{\omega+2\pi\ell}{a}\right)\right|^2<\infty .$$

Тогда, для любой функции $f\in W$, $\forall k\in Z$ выполняется неравенство

$$\| f^{(k)}\|_2^2\le \sup_{\omega}\sum_{\ell\in R}\left|\frac{\omega+2\pi\ell}{a}\right|^{2k}\left|\hat{\varphi}\left(\frac{\omega+2\pi\ell}{a}\right)\right|^2\| f\|_2^2 \qquad (3)$$

**Доказательство.** Учитывая то, что $\forall g(x)=f(ax+b)\in W$,

$$\hat{g}(\omega)=\frac{1}{a}e^{-\frac{i\omega b}{a}}\hat{f}\left(\frac{\omega}{a}\right),$$

получим

$$\| f^{(k)}\|_2^2=\frac{1}{2\pi}\left\|\left(\frac{i\omega}{a}\right)^k\hat{f}\left(\frac{\omega}{a}\right)\right\|_2^2=\frac{1}{2\pi}\int_R |i\omega|^{2k}a^{-2k}\left|\hat{f}\left(\frac{\omega}{a}\right)\right|^2 d\omega=$$

$$=\frac{1}{2\pi}\int_R \omega^{2k}a^{-2k}\left|\sum_{\gamma\in Z}c_\gamma e^{-i\gamma\omega}\hat{\varphi}\left(\frac{\omega}{a}\right)\right|^2 d\omega=$$

$$=\frac{1}{2\pi}a^{-2k}\sum_{\ell\in R}\int_0^{2\pi}|\omega+2\pi\ell|^{2k}a^2\left|\sum_{\gamma\in Z}c_\gamma e^{-i\gamma\omega}\hat{\varphi}(\omega+2\pi\ell)\right|^2 d\omega=$$

$$=\frac{1}{2\pi}a^{-2k+2}\int_0^{2\pi}|m_f(\omega)|^2\sum_{\ell\in Z}\left|\frac{\omega+2\pi\ell}{a}\right|^{2k}d\omega .$$

Заметим, что

$$\| f\|_2^2=\frac{1}{2\pi}a^2\int_0^{2\pi}|m_f(\omega)|^2 d\omega .$$

Таким образом, получим

$$\| f^{(k)}\|_2^2\le \sup_{\omega}\sum_{\ell\in R}\left|\frac{\omega+2\pi\ell}{a}\right|^{2k}\left|\hat{\varphi}\left(\frac{\omega+2\pi\ell}{a}\right)\right|^2\| f\|_2^2 .$$

Точность неравенства (3) доказывается также, как в неравенстве (1).

Теорема доказана.

Обобщим теперь неравенство (1) на многомерном случае.

Будем рассматривать подпространство $\mathrm{W}^{\mathrm{n}}$ пространства $L_2\left(R^n\right)$, порожденное сдвигами фнкции $\varphi\left(x_1,\ldots,\mathrm{x}_{\mathrm{n}}\right)$. Все свойства функций $\varphi$ и $f$ в многомерном случае остаются такими же как и в одномерном.

**Теорема 3.** Зададим функцию $\varphi\left(x_1,\ldots,\mathrm{x}_{\mathrm{n}}\right)$ такую, что $\forall \ell=\left(\ell_1,\ldots,\ell_n\right)\in R^n$, $\forall \omega=\left(\omega_1,\ldots,\omega_n\right)\in R^n$

$$\sup_{\omega_1,\ldots,\omega_n}\sum_{\ell_1,\ldots,\ell_n}\prod_{s=1}^{n}\left|\omega_s+2\pi\ell_s\right|^{2k_s}\left|\hat{\varphi}\left(\omega_1+2\pi\ell_1,\ldots,\omega_n+2\pi\ell_n\right)\right|^2<\infty .$$

Тогда, для любой функции $f\in W^{\mathrm{n}}$, $\forall k=\left(\mathrm{k}_1,\ldots,k_n\right)\in \mathrm{Z}^n$ выполняется неравенство

$$\left\|\frac{\partial^k f}{\partial x_1^{k_1}\ldots\partial x_n^{k_n}}\right\|_2^2\leq\sup_{\omega_1,\ldots,\omega_n}\sum_{\ell_1,\ldots,\ell_n}\prod_{s=1}^{n}\left|\omega_s+2\pi\ell_s\right|^{2k_s}\left|\hat{\varphi}\left(\omega_1+2\pi\ell_1,\ldots,\omega_n+2\pi\ell_n\right)\right|^2\|f\|_2^2 . \qquad (4)$$

**Доказательство.** В данном случае имеет место обобщенное равенство Парсеваля:

$$\left\|\frac{\partial^k f}{\partial x_1^{k_1}\ldots\partial x_n^{k_n}}\right\|_2^2=\left\|\omega_1^{k_1}\ldots\omega_n^{k_{n1}}\hat{f}\left(\omega\right)\right\|_2^2=\int_{R^n}\prod_{s=1}^{n}\left|\omega_s\right|^{2k_s}\left|\hat{f}\left(\omega_1,\ldots,\omega_n\right)\right|^2d\omega .$$

Так как $f\in W^{\mathrm{n}}$, то $f=\sum_{\gamma=\left(\gamma_1,\ldots,\gamma_n\right)}c_\gamma\varphi\left(x_1-\gamma_1,\ldots,x_n-\gamma_n\right)$, где $\gamma=\left(\gamma_1,\ldots,\gamma_{\mathrm{n}}\right)\in \mathrm{Z}^n$.

Оценим норму

$$\left\|\frac{\partial^k f}{\partial x_1^{k_1}\ldots\partial x_n^{k_n}}\right\|_2^2=\int_{R^n}\prod_{s=1}^{n}\left|\omega_s\right|^{2k_s}\left|\hat{f}\left(\omega_1,\ldots,\omega_n\right)\right|^2d\omega=\int_{R^n}\prod_{s=1}^{n}\left|\omega_s\right|^{2k_s}\left|\sum_{\gamma=\left(\gamma_1,\ldots,\gamma_n\right)}c_\gamma e^{-\mathrm{i}\gamma\omega}\hat{\varphi}\left(\omega\right)\right|^2d\omega=$$

$$=\int_R\ldots\int_R\prod_{s=1}^{n}\left|\omega_s\right|^{2k_s}\left|\sum_{\gamma=\left(\gamma_1,\ldots,\gamma_n\right)}c_\gamma e^{-\mathrm{i}\gamma\omega}\hat{\varphi}\left(\omega\right)\right|^2d\omega_1\ldots d\omega_n$$

Сделаем замену переменной таким образом, чтобы получить интеграл на промежутке $\left[0;2\pi\right]^n$. Тогда,

$$\left\|\frac{\partial^k f}{\partial x_1^{k_1}\ldots\partial x_n^{k_n}}\right\|_2^2=\sum_{\ell=\left(\ell_1,\ldots,\ell_n\right)_n}\int_{\left[0,2\pi\right]^n}\prod_{s=1}^{n}\left|\omega_s+2\pi\ell_s\right|^{2k_s}\left|\sum_{\gamma=\left(\gamma_1,\ldots\gamma_n\right)}c_\gamma e^{-\mathrm{i}\gamma\omega}\right|^2\left|\hat{\varphi}\left(\omega+2\pi\ell\right)\right|^2d\omega=$$

$$= \int\limits_{[0,2\pi]^n} \sum_{\ell=(\ell_1,\dots,\ell_n)} \prod_{s=1}^{n} |\omega_s + 2\pi\ell_s|^{2k_s} \left| \sum_{\gamma=(\gamma_1,\dots\gamma_n)} c_\gamma e^{-i\gamma\omega} \right|^2 |\hat{\varphi}(\omega + 2\pi\ell)|^2 d\omega,$$

где $\omega + 2\pi\ell = (\omega_1 + 2\pi\ell_1, \dots, \omega_n + 2\pi\ell_n)$

Заметим, что

$$\|f\|_2^2 = \int\limits_{[0;2\pi]^n} \sum_{\ell_1,\dots\ell_n} |\hat{\varphi}(\omega_1 + 2\pi\ell_1, \dots, \omega_n + 2\pi\ell_n)|^2 \left| \sum_{\gamma=(\gamma_1,\dots\gamma_n)} c_\gamma e^{-i\gamma\omega} \right|^2 d\omega.$$

Таким образом, получим

$$\left\| \frac{\partial^k f}{\partial x_1^{k_1} \dots \partial x_n^{k_n}} \right\|_2^2 = \int\limits_{[0,2\pi]^n} \left| \sum_{\gamma=(\gamma_1,\dots\gamma_n)} c_\gamma e^{-i\gamma\omega} \right|^2 \sum_{\ell=(\ell_1,\dots,\ell_n)} \prod_{s=1}^{n} |\omega_s + 2\pi\ell_s|^{2k_s} |\hat{\varphi}(\omega + 2\pi\ell)|^2 d\omega \le$$

$$\le \sup_{\omega_1,\dots\omega_n} \sum_{\ell_1,\dots,\ell_n} \prod_{s=1}^{n} |\omega_s + 2\pi\ell_s|^{2k_s} |\hat{\varphi}(\omega_1 + 2\pi\ell_1, \dots, \omega_n + 2\pi\ell_n)|^2 \|f\|_2^2 =$$

$$= \int\limits_{[0,2\pi]^n} \left| \sum_{\gamma=(\gamma_1,\dots\gamma_n)} c_\gamma e^{-i\gamma\omega} \right|^2 \sum_{\ell=(\ell_1,\dots,\ell_n)} \prod_{s=1}^{n} |\omega_s + 2\pi\ell_s|^{2k_s} |\hat{\varphi}(\omega + 2\pi\ell)|^2 d\omega,$$

Покажем теперь, что константа в правой части неравенства (4) неулучшаема..

Выберем семейство функций $f_m \in W^n$ таких, что

$$\frac{\left\| \frac{\partial^k f_m}{\partial x_1^{k_1} \dots \partial x_n^{k_n}} \right\|_2^2}{\|f_m\|_2^2} \xrightarrow[n\to\infty]{} \sup_{\omega_1,\dots\omega_n} \sum_{\ell_1,\dots,\ell_n} \prod_{s=1}^{n} |\omega_s + 2\pi\ell_s|^{2k_s} |\hat{\varphi}(\omega_1 + 2\pi\ell_1, \dots, \omega_n + 2\pi\ell_n)|^2$$

Определим последовательность $f_m$ следующим образом:

положим

$$\hat{f}_m(\omega) = u_m(\omega)\hat{\varphi}(\omega),$$

где

$$u_m(\omega) = \sqrt{\tilde{\Phi}_m(\omega)} := \sqrt{\prod_{s=1}^{n} \Phi_m(\omega_s)}.$$

Рассмотрим

$$\frac{\left\| \frac{\partial^k f_m}{\partial x_1^{k_1} \dots \partial x_n^{k_n}} \right\|_2^2}{\|f_m\|_2^2} = \frac{\left\| \prod_{s=1}^{n} \omega_s^{k_s} \hat{f}_m(\omega) \right\|_2^2}{\|\hat{f}_m\|_2^2} = \frac{\left\| \prod_{s=1}^{n} \omega_s^{k_s} \sqrt{\tilde{\Phi}_n(\omega)} \hat{\varphi}(\omega) \right\|_2^2}{\left\| \sqrt{\tilde{\Phi}_n(\omega)} \hat{\varphi}(\omega) \right\|_2^2} =$$

$$=\frac{\int\limits_{R^n}\prod\limits_{s=1}^{n}\omega_s^{k_s}\tilde{\Phi}_n(\omega)\left|\hat{\varphi}(\omega)\right|^2 d\omega}{\int\limits_{R^n}\tilde{\Phi}_n(\omega)\left|\hat{\varphi}(\omega)\right|^2 d\omega}=\frac{\sum\limits_{\ell}\int\limits_{[0,2\pi]^n}\prod\limits_{s=1}^{n}\left|\omega_s+2\pi\ell_s\right|^{k_s}\tilde{\Phi}_n(\omega)\left|\hat{\varphi}(\omega+2\pi\ell)\right|^2 d\omega}{\sum\limits_{\ell}\int\limits_{[0,2\pi]^n}\tilde{\Phi}_n(\omega)\left|\hat{\varphi}(\omega+2\pi\ell)\right|^2 d\omega}=$$

$$=\frac{\int\limits_{[0,2\pi]^n}\sum\limits_{\ell}\prod\limits_{s=1}^{n}\left|\omega_s+2\pi\ell_s\right|^{k_s}\tilde{\Phi}_n(\omega)\left|\hat{\varphi}(\omega+2\pi\ell)\right|^2 d\omega}{\int\limits_{[0,2\pi]^n}\sum\limits_{\ell}\tilde{\Phi}_n(\omega)\left|\hat{\varphi}(\omega+2\pi\ell)\right|^2 d\omega}.$$

Учитывая свойства (2) ядра Фейера $\Phi_n(\omega)$, получим:

$$\frac{\left\|\frac{\partial^k f_m}{\partial x_1^{k_1}\ldots\partial x_n^{k_n}}\right\|_2^2}{\|f_m\|_2^2}=\frac{\int\limits_{[0,2\pi]^n}\sum\limits_{\ell}\prod\limits_{s=1}^{n}\left|\omega_s+2\pi\ell_s\right|^{k_s}\tilde{\Phi}_n(\omega)\left|\hat{\varphi}(\omega+2\pi\ell)\right|^2 d\omega}{\int\limits_{[0,2\pi]^n}\sum\limits_{\ell}\tilde{\Phi}_n(\omega)\left|\hat{\varphi}(\omega+2\pi\ell)\right|^2 d\omega}\xrightarrow[n\to\infty]{}$$

$$\xrightarrow[n\to\infty]{}\sup_{\omega}\sum_{\ell}\prod_{s=1}^{n}\left|\omega_s+2\pi\ell_s\right|^{2k_s}\left|\hat{\varphi}(\omega+2\pi\ell)\right|^2$$

Теорема доказана.

**Литература**